\documentclass{icmart}
\usepackage{verbatim}
\RequirePackage[colorlinks,citecolor=blue,urlcolor=blue]{hyperref}
\usepackage{breakurl}
\usepackage{graphicx,psfrag}
\usepackage{mathtools} 
\usepackage{wasysym}   

\def\mik{1}
\newcommand\cpsfrag[2]{\ifnum\mik=1\psfrag{#1}{#2}\fi}

\newtheorem{theorem}{Theorem}
\newtheorem{conjecture}[theorem]{Conjecture}

\newtheorem{corollary}[theorem]{Corollary}

\numberwithin{equation}{section}
\numberwithin{theorem}{section}
\numberwithin{figure}{section}
\numberwithin{table}{section}

\newcounter{mycount}

\newenvironment{numlist}{\begin{list}{\arabic{mycount}.}%
   {\usecounter{mycount}\labelwidth=1cm\itemsep 0pt}}{\end{list}}
\newenvironment{letlist}{\begin{list}{\rm(\alph{mycount})}%
   {\usecounter{mycount}\labelwidth=1cm\itemsep 0pt}}{\end{list}}
\newenvironment{problist}{\begin{list}{{\alph{mycount}}}%
   {\usecounter{mycount}\leftmargin=0pt\itemindent=1cm\labelwidth=1cm\itemsep 0pt}}{\end{list}}

\newcommand\bxp{box-crossing property}
\newcommand\stt{star--triangle transformation}
\newcommand\rad{{\text{\rm rad}}}
\newcommand\dc{d_{\text{\rm c}}}
\newcommand\rc{random-cluster}

\newcommand\oo{\infty}
\newcommand\De{\Delta}

\newcommand\TT{{\mathbb T}}
\newcommand\EE{{\mathbb E}}
\newcommand\sL{{\mathcal L}}

\newcommand\sG{{\mathcal G}}

\newcommand\ZZ{{\mathbb Z}}

\newcommand\RR{{\mathbb R}}

\newcommand\om{\omega}

\renewcommand\a{\alpha}

\newcommand\eps{\epsilon}
\renewcommand\b{\beta}

\newcommand\resp{respectively}

\newcommand\de{\delta}

\newcommand\q{\quad}

\newcommand\pd{\partial}

\newcommand\qq{\qquad}
\newcommand\Om{\Omega}
\newcommand\PP{{\mathbb P}}
\newcommand\lra{\leftrightarrow}
\newcommand\pc{p_{\text{\rm c}}}

\newcommand\bp{\mathbf{p}}

\newcommand\fpq{\phi_{p,q}}

\newcommand\fgq{\phi_{G,q}}
\newcommand\fgbq{\phi_{G,q,\b}}

\newcommand\fpqb{\fpq^b}

\newcommand\pcb{\pc^b}

\contact[g.r.grimmett@statslab.cam.ac.uk]{Geoffrey R.\ Grimmett}

\title{Criticality, universality, and isoradiality}
\author{Geoffrey R.\ Grimmett}

\begin{document}
\begin{abstract}
Critical points and singularities are encountered in the study of critical 
phenomena in probability and physics. We present recent results concerning
the values of such critical points 
and the nature of the singularities for two prominent probabilistic models,
namely percolation and the more general random-cluster model.
The main topic is the statement and proof of the criticality and universality
of the canonical measure of bond percolation on
isoradial graphs (due to the author and Ioan Manolescu). The key technique used in this work is the \stt,
known also as the Yang--Baxter equation.  The second topic reported here is
the identification of the critical point of the \rc\ model on the square lattice
(due to Beffara and Duminil-Copin), and of
the criticality of the canonical measure
of the \rc\ model with $q \ge 4$ on periodic isoradial graphs (by the same authors with
Smirnov). The proof of universality for percolation is expected
to extend to the \rc\ model on isoradial graphs. 
\end{abstract}


\begin{keywords}
Percolation, random-cluster model, Ising/Potts models,
critical point, universality, 
isoradial graph, critical exponent, \stt, Yang--Baxter equation.
\end{keywords}

\begin{classification}
Primary 60K35; Secondary 82B30.
\end{classification}

\maketitle

\section{Introduction}\label{sec:intro}

One of the most provocative and elusive problems in the mathematics of critical phenomena
is the issue of \emph{universality}. Disordered physical systems manifest phase transitions,
the nature of which is believed to be independent of the local structure of space.  
Very little about universality is known rigorously for systems below their upper critical
dimension. 
It is frequently said that \lq\lq renormalization" is the key to  universality, but rigorous
applications of renormalization in the context of universality are rare.

There has been serious recent progress in the \lq\lq exactly solvable" 
setting of the two-dimensional Ising model, and a handful of special cases for other models.
Our principal purpose here is to outline recent progress
concerning the identification of critical surfaces and the issue of universality for
bond percolation and the \rc\ model on isoradial graphs, with emphasis on the
general method, the current limitations, and the open problems.

For bond percolation on an extensive family of isoradial graphs, the canonical process,
in which the \stt\ is in harmony with the geometry, is shown to be critical.
Furthermore, universality has been proved for this class of systems, at least for
the critical exponents \emph{at the critical surface}. These results, found in 
recent papers by the author and Manolescu, \cite{GM1,GM2,GM3},
vastly extend earlier calculations of critical values for the square lattice etc, with
the added ingredient of universality. Note that, to date, we are able to prove only \emph{conditional}
universality: if a certain exponent exists for at least one isoradial graph, 
then a family of exponents exist for an extensive
collection of isoradial graphs, and they are universal across this collection.

The picture for the general \rc\ model is more restrained, but significant
progress has been achieved on the identification of critical points. 
The longstanding conjecture for the critical value of the square lattice
has been proved by Beffara and Duminil-Copin \cite{Beffara_Duminil}, using a development
of classical tools. Jointly with Smirnov \cite{BDS2}, the same authors have
used Smirnov's parafermionic observable in the first-order
setting of $q\ge 4$ to identify the critical surface of a periodic isoradial graph.
It is conjectured that the methods of \cite{GM3} may be extended to
obtain universality for the \rc\ model on isoradial graphs.

The results reported in this survey are closely related to certain famous `exact results' 
in the physics literature. Prominent in the latter regard is the book of Baxter \cite{Baxter_book},
from whose preface we quote selectively as follows:
\begin{quotation}
\lq\lq $\dots$ the phrase `exactly solved' has been chosen with care. It is not necessarily
the same as `rigorously solved'. $\dots$ 
There is of course still much to be done."
\end{quotation}

Percolation is summarized in Section \ref{sec:perc}, and isoradial graphs in
Section \ref{sec:iso}. Progress with criticality and universality for percolation are described in
Section \ref{sec:crit+univ}. Section \ref{sec:rcm} is devoted to recent progress
with critical surfaces of \rc\ models on isoradial graphs, and open problems 
for percolation and the \rc\ model may be found
in Sections \ref{sec:open-p} and \ref{sec:open-rcm}. 

\section{Percolation}\label{sec:perc}

\subsection{Background}

Percolation is the fundamental stochastic model for spatial disorder.
Since its introduction by Broadbent and Hammersley in 1957,
it has emerged as a key topic in probability theory, with connections and impact across
all areas of applied science in which disorder meets geometry. 
It is in addition a source of beautiful and apparently difficult
mathematical problems, the solutions to which often require the development
of new tools with broader applications.

Here is the percolation process in its basic form.
Let $G=(V,E)$ be an infinite, connected graph, typically
a crystalline lattice such as the $d$-dimensional hypercubic lattice. We are provided
with a coin that shows heads with some fixed probability $p$.
For each edge $e$ of $G$, we flip the coin, and we designate $e$  \emph{open} if
heads shows, and \emph{closed} otherwise.  The open edges are considered open to the
passage of material such as liquid, disease, or rumour.\footnote{This is the process
known as \emph{bond} percolation.  Later we shall refer to \emph{site} percolation, in which
the vertices (rather than the edges) receive random states.}

Liquid is supplied at a \emph{source} vertex $s$, and it flows 
along the open edges and is blocked by the closed edges.
The basic problem is to determine the geometrical properties
(such as size, shape, and so on) of the region $C_s$ that is wetted by the liquid.
More generally, one is interested in the geometry of
the connected subgraphs of $G$ induced by the set of open edges. The components
of this graph are called the \emph{open clusters}.

Broadbent and Hammersley proved in \cite{BrH,H57a,H59} 
that there exists a \emph{critical probability}
$\pc=\pc(G)$ such that: every open cluster is bounded if $p< \pc$, and
some  open cluster is unbounded if $p> \pc$.  There are
two \emph{phases}:
the \emph{subcritical phase} when $p< \pc$ and the \emph{supercritical phase}
when $p> \pc$. The singularity that occurs when $p$ is near or equal to $\pc$
has attracted a great deal of attention from mathematicians and physicists,
and many of the principal problems remain unsolved even after several decades
of study.  See \cite{G99,Grimmett_Graphs}
for general accounts of the theory of percolation.

Percolation is one of a large family of models of 
classical and quantum statistical physics that manifest phase transitions,
and its theory is near the heart of the extensive scientific project to
understand phase transitions and critical phenomena. Key aspects of
its special position in the general theory include: 
(i) its deceptively simple formulation as a probabilistic model, (ii) its
use as a comparator for more complicated systems, and (iii)
its role in the development of new methodology. 

One concrete connection between percolation and models for ferromagnetism
is its membership of the one-parameter family of so-called random-cluster models.
That is, percolation is the $q=1$ \rc\ model. The $q=2$ \rc\ model corresponds
to the Ising  model, and the $q=3,4,\dots$ \rc\ models to the $q$-state Potts models.
The $q \downarrow 0$ limit is connected to electrical networks,
uniform spanning trees, and uniform connected subgraphs.
The \emph{geometry} of the \rc\ model corresponds
to the \emph{correlation} structure of the Ising/Potts models,
and thus its critical point $\pc$ may be expressed in terms of
the critical temperature of the latter systems. See \cite{G-rcm,Werner_SMF}
for a general account of the \rc\ model.

The theory of percolation is extensive and influential. Not only is percolation a
benchmark model for studying random spatial processes in general, but also it has been,
and continues to be, a source of intriguing and beautiful open problems.
Percolation in two dimensions has been especially prominent in the last decade by virtue
of its connections to conformal invariance and conformal field theory. Interested readers
are referred to the papers \cite{Cardy,G-rev,Sch06,Smirnov, Smi07, Sun11, WW_park_city} and the books
\cite{BolRio, G99,Grimmett_Graphs}.

\subsection{Formalities}\label{sec:form}

For $x,y\in V$, we write $x \lra y$ if there exists an open path
joining $x$ and $y$. The {\it open cluster\/} at the vertex $x$ is the
set $C_x =\{y: x\lra y\}$ of all
vertices reached along open paths from $x$,
and we write $C=C_0$ where $0$ is a fixed vertex
called the \emph{origin}.
Write $\PP_p$ for the relevant product probability measure,
and $\EE_p$ for expectation with respect to $\PP_p$.

The {\it percolation probability\/} is the function $\theta(p)$ given by
$$
\theta(p) = \PP_p(|C|=\oo),
$$
and the \emph{critical probability} is defined by
\begin{equation}\label{defcritprob}
\pc = \pc(G) = \sup\{p: \theta(p)=0\}.
\end{equation}
It is elementary that $\theta$ is a
non-decreasing function, and therefore,
$$
\theta(p)\, \begin{cases} = 0 &\text{if } p<\pc,\\
>0 &\text{if } p>\pc.
\end{cases}
$$
It is not hard to
see, by the Harris--FKG inequality,  that the value $\pc(G)$ does not depend on the choice of origin.

Let $d \ge 2$, and let $\sL$ be a $d$-dimensional lattice.
It is a fundamental fact that $0<\pc(\sL)<1$, but it is unproven in general that
no infinite open cluster exists when $p=\pc$.

\begin{conjecture}\label{thetapc0}
For any lattice $\sL$ in $d \ge 2$ dimensions, we have that
$\theta(\pc) = 0$.
\end{conjecture}

The claim of the conjecture is known to be valid for certain lattices when $d=2$
and for large $d$, currently $d \ge 15$. 
This conjecture has been the `next open problem' since the intensive study of the late 1980s.

Whereas the above process is defined in terms
of a single parameter $p$,
we are concerned here with the richer multi-parameter setting in which
an edge $e$ is designated open with some probability $p_e$.
In such a case, the critical probability $\pc$ is replaced by a so-called `critical surface'.

\subsection{Critical exponents and universality}\label{sec:pt}

A great deal of effort has been directed towards
understanding the nature of the percolation
phase transition. The picture is now fairly clear for one specific model
in two dimensions (site percolation on the triangular lattice),
owing to
the very significant progress in recent years linking critical percolation
to the Schramm--L\"owner curve SLE$_6$. There remain however
substantial difficulties to be overcome even when $d=2$,
associated largely with the extension of such results to general two-dimensional systems.
The case of large $d$ (currently, $d \ge 15$) is also
well understood, through work based on the so-called `lace expansion' (see \cite{BDGS}).
Many problems remain open in the prominent case $d=3$.

Let $\sL$ be a $d$-dimensional lattice.
The nature of the percolation singularity on $\sL$  is
expected to  share general
features with phase transitions of other models of statistical mechanics.
These features are sometimes referred to as `scaling theory' and they
relate to the `critical exponents' occurring in the power-law
singularities (see \cite[Chap.\ 9]{G99}).
There are two sets of critical exponents, arising firstly in the limit as
$p\to\pc$, and secondly in the limit over increasing
spatial scales when $p=\pc$. The definitions of the critical exponents
are found  in Table~\ref{Tab-ce} (taken from \cite{G99}).

\begin{table}[t]
\centering\small
\tabcolsep=5pt

\begin{tabular}{|c|c|c|c|}
\hline
\multicolumn{2}{|c|}{}&&\\
\multicolumn{2}{|c|}{\emph{Function}}  & \emph{Behaviour} & \emph{Exp.} \\
\multicolumn{2}{|c|}{}&&\\
\hline
&&&\\
\raisebox{1ex}{percolation}& $\theta (p)=\PP_p(|C|=\infty )$ & $\theta (p)\approx (p-\pc )^\beta$ & $\beta$ \\
\raisebox{.8ex}{probability}&&&\\
&&&\\
\raisebox{1ex}{truncated}& $\chi^{\text f}(p)=\EE_p(|C|1_{|C|<\infty})$& $\chi^{\text f}(p)\approx |p-\pc |^{-\gamma}$ & $\gamma$\\
\raisebox{.8ex}{mean cluster-size} &&& \\
&&&\\
\raisebox{1ex}{number of}& $\kappa (p)=\EE_p(|C|^{-1})$ & $\kappa'''(p)\approx |p-\pc |^{-1-\alpha}$ & $\alpha$\\
\raisebox{.8ex}{clusters per vertex}&&&\\
&&&\\
cluster moments& $\chi_k^{\text f} (p)=\EE_p(|C|^k1_{|C|<\infty} )$ & $\displaystyle\frac{\chi_{k+1}^{\text f}
     (p)}{\chi_k^{\text f}(p)}\approx |p-\pc |^{-\De}$ & $\De$\\
&&&\\
correlation length& $\xi (p)$ & $\xi (p)\approx |p-\pc |^{-\nu}$ & $\nu$\\
&&&\\
\hline
\multicolumn{2}{|c|}{}&&\\
\multicolumn{2}{|c|}{cluster volume}& $\PP_{\pc} (|C|=n)\approx n^{-1-1/\de}$ & $\de$\\
\multicolumn{2}{|c|}{}&&\\
\multicolumn{2}{|c|}{cluster radius} & $\PP_{\pc}\bigl(\rad (C)=n\bigr)\approx n^{-1-1/\rho}$ & $\rho$\\
\multicolumn{2}{|c|}{}&&\\
\multicolumn{2}{|c|}{connectivity function} & $\PP_{\pc}(0\lra x)\approx   \| x\|^{2-d-\eta}$ & $\eta$\\
\multicolumn{2}{|c|}{}&&\\
\hline
\end{tabular}
\caption{Eight functions and their critical exponents. The first five exponents arise
in the limit as $p \to \pc$, and the remaining three as $n\to\oo$ with $p=\pc$.
See \cite[p.\ 127]{G99} for a definition
 of the correlation length $\xi(p)$.}
\label{Tab-ce}
\end{table}

The notation of Table~\ref{Tab-ce} is as follows.
We write $f(x) \approx g(x)$ as $x \to x_0 \in[0,\oo]$ if
$\log f(x)/\log g(x) \to 1$.
 The \emph{radius} of the open cluster $C$ at the origin $x$ is defined by
$$
\rad (C)
=\sup\{\|y\|: x \lra y\},
$$
where
$$
\|y\|= \sup_i |y_i|, \qquad y=(y_1,y_2,\dots,y_d)\in\RR^d,
$$
is the supremum ($L^\oo$) norm on $\RR^d$.
The limit as $p\to\pc$ should be interpreted
in a manner appropriate for the function in question
(for example, as $p \downarrow \pc$ for $\theta(p)$,
but as $p\to\pc$ for $\kappa(p)$).
The \emph{indicator function} of an event $A$ is denoted $1_A$.

Eight critical exponents are listed in Table~\ref{Tab-ce},
denoted $\alpha$, $\beta$, $\gamma$, $\de$, $\nu$, $\eta$, $\rho$, $\De$,
but there is no general proof of the
existence of any of these exponents for arbitrary $d\ge 2$.
Such critical exponents may be defined for phase
transitions in a large family of physical systems. The exponents
are not believed to be independent variables, but rather to satisfy the so-called {\it scaling relations\/}
\begin{gather*}
2-\alpha=\gamma +2\beta =\beta (\de +1),\\
\De=\de\beta ,\quad
\gamma=\nu (2-\eta ),
\end{gather*}
and, when $d$ is not too large, the {\it hyperscaling relations\/}\label{ind-hypers}
\begin{align*}
d\rho=\de +1,\quad
2-\alpha=d\nu .
\end{align*}

More generally, a `scaling relation' is any equation involving critical
exponents which is believed to be `universally' valid.
The {\it upper critical dimension\/} is the largest value $\dc$
such that the hyperscaling relations hold for $d\leq\dc$ and not otherwise. 
It is believed
that $\dc =6$ for percolation.
There is no general proof of the validity of the scaling and hyperscaling
relations for percolation, although quite a lot is known when either $d=2$ or $d$ is large.
The case of large $d$ is studied via the lace expansion, and this
is expected to be valid for $d > 6$.

We note some further points in the context of percolation.
\begin{letlist}
\item
{\it Universality\/}. The numerical values of critical
exponents are believed to depend only on the value of $d$, and to be independent of the
choice of lattice, and of the type of percolation under study. 

\item
{\it Two dimensions\/}. When $d=2$, it is believed that
$$
\a =-\tfrac23,\ \b =\tfrac5{36},\ \gamma =\tfrac{43}{18},\ \de =\tfrac{91}5
,\dots .
$$
These values (other than $\a$) have been proved  (essentially only) in the special case of site
percolation on the triangular lattice, see \cite{LSW01,Smirnov-Werner}.

\item
{\it Large dimensions\/}. When $d$ is sufficiently large (in fact,
$d\geq\dc$) it is believed that the critical exponents are the same as
those for percolation on a tree (the `mean-field model'), namely $\de =2$,
$\gamma = 1$, $\nu=\frac12$,
$\rho =\frac12$, and so on. Using the
first hyperscaling relation, this is consistent with the contention that $\dc =6$.
Several such statements are known to hold for $d\ge 15$, see \cite{Fitz,HS, HS94, KN}.
\end{letlist}

Open challenges include the following:
\begin{numlist}
\item prove the existence of critical exponents for general lattices,
\item prove some version of universality,
\item prove the scaling and hyperscaling relations in general dimensions,
\item calculate the critical exponents for general models in two dimensions,
\item prove the mean-field values of critical exponents when $d \ge 6$.
\end{numlist}
Progress towards these goals has been substantial but patchy. As noted above,
for sufficiently large $d$, the lace expansion has enabled proofs
of exact values for many exponents, for a restricted class of lattices. 
There has been remarkable progress
in recent years when $d=2$, inspired
largely by work of Cardy \cite{Cardy} and Schramm \cite{Sch00}, enacted by Smirnov \cite{Smirnov},
and confirmed by the programme pursued by Lawler, Schramm, Werner,
Camia, Newman, Sheffield and others to understand SLE curves and conformal ensembles.

In this paper, we concentrate on recent progress concerning isoradial embeddings 
of planar graphs, and particularly the identification of their critical surfaces 
and the issue of universality.

\section{Isoradial graphs}\label{sec:iso}

Let $G$ be an infinite, planar graph embedded in $\RR^2$ in such a way that edges intersect only
at vertices. For simplicity, we assume that the embedding has only bounded faces.
The graph $G$ is called \emph{isoradial} if (i) every face has
a circumcircle which passes through every vertex of the face,  (ii)
the centre of each circumcircle lies in the interior
of the corresponding face, and (iii) all such circumcircles have the same radius.
We may assume by re-scaling that the common radius is $1$.  

\begin{figure}[t]
 \centering
    \includegraphics[width=0.47\textwidth]{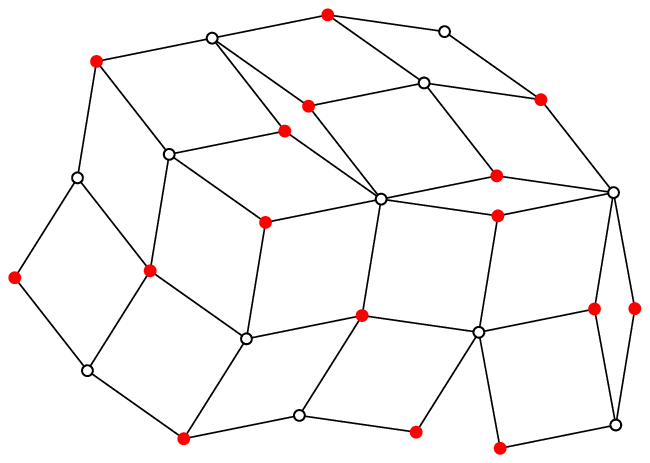}
\q\includegraphics[width=0.47\textwidth]{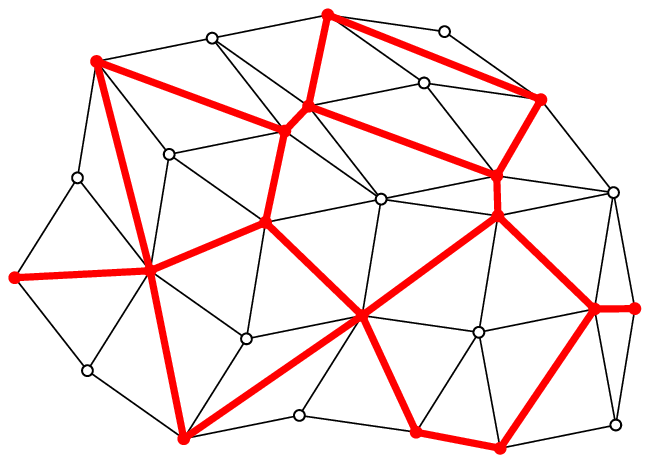}  
\caption{On the left is part of a rhombic tiling of the plane. Since all cycles
have even length, this is a bipartite graph, with vertex-sets coloured red and white. 
The graph on the right is obtained by joining pairs of red vertices
across faces. Each red face of the latter graph contains a unique white vertex,
and this is the centre of the circumcircle of that face. 
Joining the white vertices, instead, yields
another isoradial graph that is dual to the first.}
  \label{fig:rho}
\end{figure}

\begin{figure}[t]
 \centering
    \includegraphics[width=0.4\textwidth]{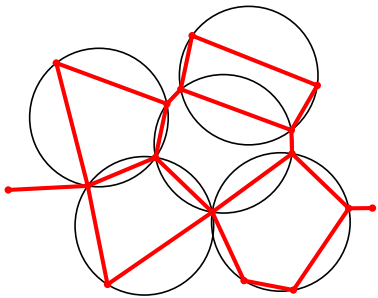}
\caption{An illustration of the isoradiality of the red graph of Figure
\ref{fig:rho}.}
  \label{fig:rho2}
\end{figure}

The family of isoradial graphs is in two-to-one correspondence with the family
of tilings of the plane with rhombi of side-length $1$, in the following sense.
Consider a rhombic tiling of the plane, as in Figure \ref{fig:rho}.  The tiling, when viewed
as a graph, is bipartite with vertex-sets coloured red and white, say. Fix a colour
and join any two vertices of that colour whenever they are the opposite vertices of
a rhombus. The resulting graph $G$ is isoradial.  If the other colour is chosen, the resulting
graph is the (isoradial) dual of $G$.
 This is illustrated in Figures
\ref{fig:rho} and \ref{fig:rho2}.
Conversely, given an isoradial graph $G$, the corresponding rhombic tiling is
obtained by augmenting its vertex-set by the circumcentres of the faces, and
each circumcentre is joined to the vertices of the enclosing face.

Isoradial graphs were introduced by Duffin \cite{Duff}, and are related
to the so-called $Z$-invariant graphs of Baxter \cite{Bax78}. They were named thus 
by Kenyon, whose expository paper \cite{Ken02} proposes the 
connection between percolation and isoradiality (and much more).
Isoradial graphs have two important properties, the first of
which is their connection to preholomorphic functions. This was discovered
by Duffin, and is summarized by Smirnov  \cite{Smi-icm} and 
developed further in the context of probability by Chelkak and Smirnov \cite{Chelkak-Smirnov3}. 
This property is key to the work on the \rc\ model on isoradial graphs reviewed
in Section \ref{sec:rcm}. A recent review of connections between isoradiality and 
aspects of statistical mechanics may be found in \cite{BdT12}.

The second property of isoradial graphs
is of special relevance in the current work, namely that
they provide the `right' setting for the \stt. This is explained next.

\begin{figure}[t]
 \begin{center}\normalsize
	\cpsfrag{th}{$\theta_e$}
	\cpsfrag{e}{$e$}
\includegraphics[width=0.6\textwidth]{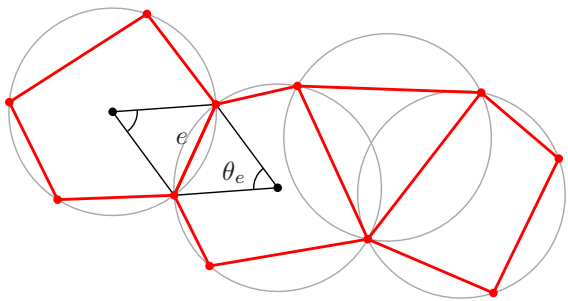}
\end{center}
\caption{The edge $e$ is the diagonal of some rhombus, with opposite angle
$\theta_e$ as illustrated.}
  \label{fig:isomap}
\end{figure}

Consider an inhomogeneous bond percolation process on the isoradial graph $G$,
whose edge-probabilities $p_e$ are given as follows in terms of the graph-embedding. Each edge $e$
of $G$ is the diagonal of a unique rhombus in the corresponding rhombic tiling of
the plane, and its parameter $p_e$ is given in terms
of the geometry of this rhombus.  With $\theta_e$ the opposite angle of the rhombus,
as illustrated in Figure \ref{fig:isomap}, let $p_e \in (0,1)$ satisfy
\begin{equation}\label{pdef}
\frac{p_e}{1-p_e} = \frac{\sin(\frac13[\pi-\theta_e])}{\sin(\frac13\theta_e)}.
\end{equation}
We consider inhomogeneous bond percolation on $G$ in which each edge $e$ is designated open
with probability $p_e$, and we refer to this as the
\emph{canonical percolation process} on $G$, with associated probability
measure $\PP_G$. The special property of the vector
$\bp=(p_e: e\in E)$ is explained in Section \ref{sec:stt}.

In a beautiful series of  papers \cite{deB1,deB2,deB3}, de Bruijn introduced the
geometrical construct of `ribbons' or `train tracks' via which he was able
to build a theory of rhombic tilings. Consider a tiling $\TT$ of the plane in which each tile
is convex with four sides. We pursue a walk on the faces of $\TT$ according to the following rules. 
The walk starts in some given tile, and crosses some edge to a neighbouring tile.
It next traverses the opposite edge of this tile, and so on. The walk may
be extended backwards according to the same rule, and a doubly-infinite walk ensues.
Such a walk is called a \emph{ribbon} or \emph{track}. De Bruijn pointed out
that, if $\TT$ is a rhombic tiling, then no walk intersects itself, and two walks may 
intersect once but not twice. This property turns out to be both necessary and sufficient
for a track system to be homeomorphic to that of a \emph{rhombic} tiling (see \cite{KenS}). 

\begin{figure}[t]
 \begin{center}
\includegraphics[width=0.5\textwidth]{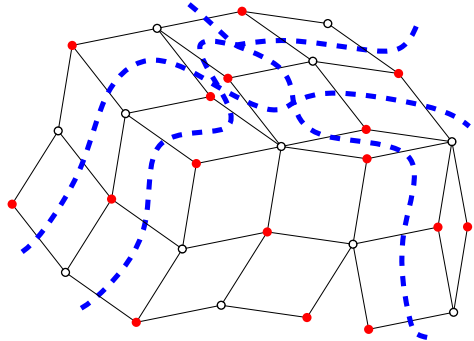}
\end{center}
\caption{An illustration of the track system of the rhombic tiling of Figure
\ref{fig:rho}.}
  \label{fig:track}
\end{figure}

We impose two restrictions on the isoradial graphs under study.
Firstly, we  say that an isoradial graph $G=(V,E)$ satisfies the \emph{bounded-angles property}
(BAP) if there exists $\eps>0$ such that
$$
\eps < \theta_e < \pi-\eps\qq\text{for all } e \in E,
$$
where $\theta_e$ is as in Figure \ref{fig:isomap}.
This amounts to the condition that the rhombi in the corresponding tiling are not
`too flat'. We say that $G$ has the \emph{square-grid property} (SGP) if
its track system, viewed as a graph, contains a square grid such that those
tracks not in the grid have boundedly many intersections with the grid within any bounded 
region (see \cite[Sect.\ 4.2]{GM3} for a more careful statement of this property).

An isoradial graph may be viewed as both a graph and a planar embedding of a graph. 
Of the many examples of isoradial graphs, we mention first
the conventional embeddings of the square, triangular, and hexagonal lattices. These are
symmetric embeddings, and the edges have the same $p$-value. There are also
non-symmetric isoradial embeddings of the same lattices, 
and indeed embeddings with no non-trivial symmetry,
for which the corresponding percolation measures are `highly inhomogeneous'.

\begin{figure}
\begin{center}
\includegraphics[width=0.47\textwidth]{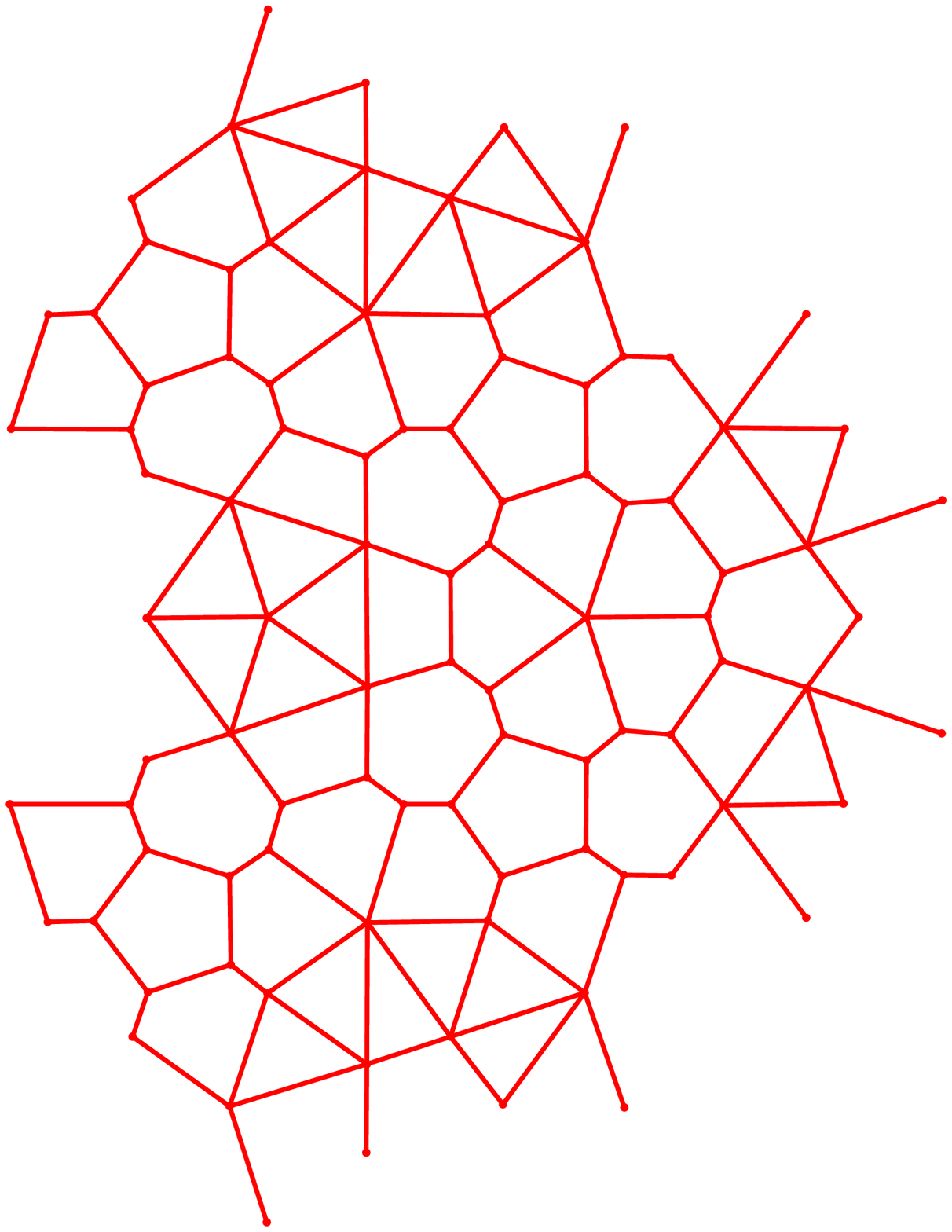}
\quad \includegraphics[width=0.47\textwidth]{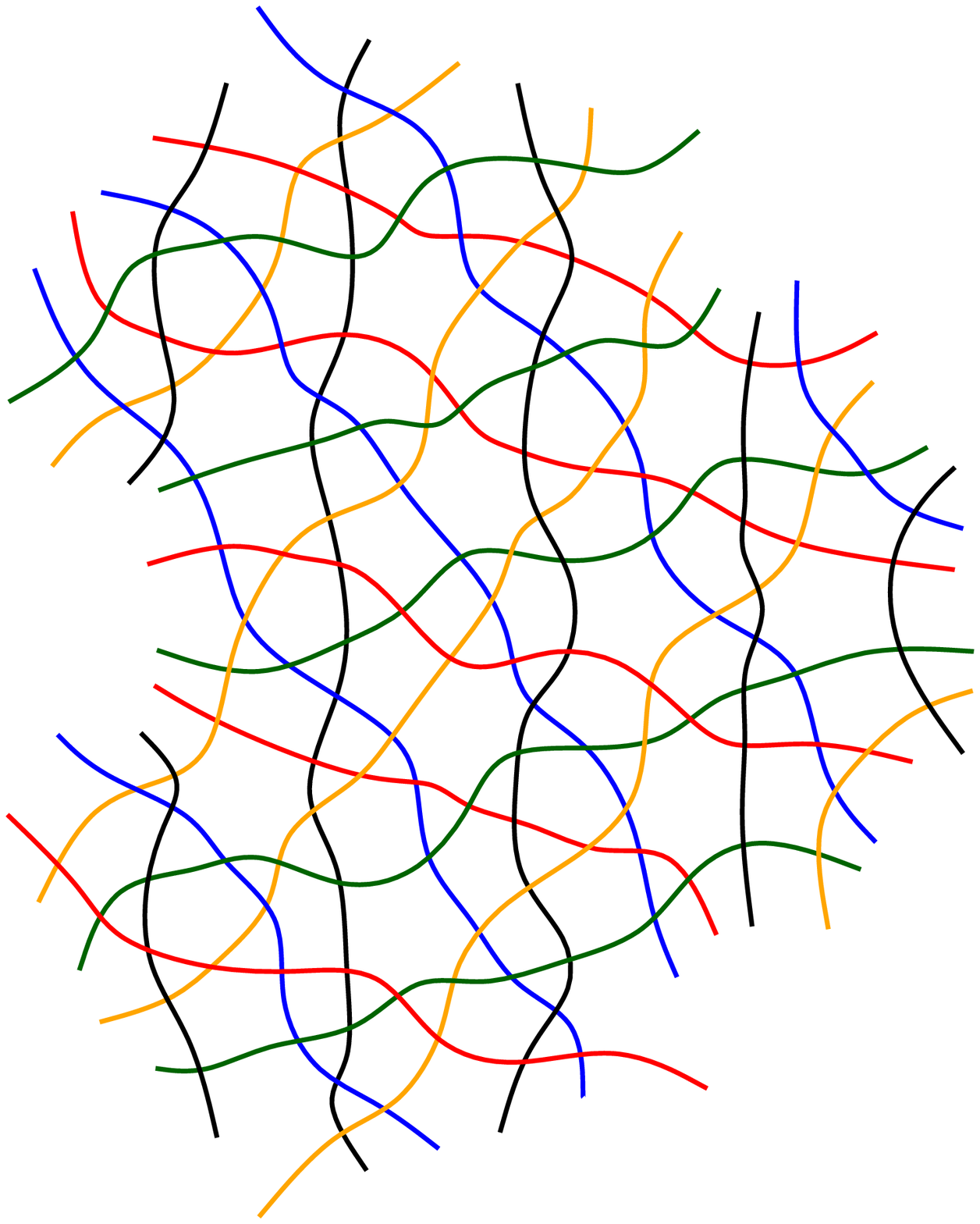}
\end{center}
\caption{On the left, an isoradial graph obtained from part of the Penrose rhombic tiling.
On the right, the associated track system comprises a pentagrid: five sets of
non-intersecting doubly-infinite lines.}
\label{fig:pen}
\end{figure}

The isoradial family is much richer than the above examples might indicate, and includes
graphs obtained from aperiodic tilings including the classic Penrose tiling \cite{Pen74,Pen78},
illustrated in Figure \ref{fig:pen}. All isoradial graphs mentioned above satisfy the SGP,
and also the BAP so long as the associated tiling comprises rhombi 
with flatness uniformly bounded from $0$.

\section{Criticality and universality for percolation}\label{sec:crit+univ}

\subsection{Two main theorems}

The first main theorem of \cite{GM3} is the identification of
the criticality
of the canonical percolation measure $\PP_G$ on an isoradial graph $G$. The
second is the universality of $\PP_G$ across an extensive family of isoradial graphs $G$. 

In order to state the criticality theorem, we introduce notation that is appropriate for a perturbation of the 
canonical measure $\PP_G$, and we borrow that of \cite{BDS2}. For
$e\in E$ and $\b\in(0,\oo)$, let $p_e(\b)$ satisfy
\begin{equation}\label{pdef2}
\frac{p_e(\b)}{1-p_e(\b)} = \b\frac{\sin(\frac13[\pi-\theta_e])}{\sin(\frac13\theta_e)},
\end{equation}
and write $\PP_{G,\b}$ for the corresponding product measure on $G$.
Thus $\PP_{G,1}=\PP_G$.

\begin{theorem}[Criticality \cite{GM3}]\label{thm:1}
Let $G=(V,E)$ be an isoradial graph with the bounded-angles property and the square-grid property. 
The canonical percolation measure $\PP_G$ is critical in
that
  \begin{letlist}
  \item
    there exist $a,b,c,d > 0$ such that
    \begin{align*}
      ak^{-b} \le  \PP_{G} \bigl(\rad (C_v) \geq k\bigr) \leq ck^{-d}, 
      \qquad k \ge 1, \q v \in V,
    \end{align*}
  \item
    there exists, $\PP_G$-a.s.,  no infinite open cluster, 
  \item
    for $\b<1$, there exist $f,g>0$ such that
    $$
    \PP_{G,\b}(|C_v|\ge k) \le f e^{-gk}, \qquad k \ge 0,\q v \in V,
    $$
   \item 
    for $\b>1$, there exists, $\PP_{G,\b}$-a.s., a unique infinite open cluster. 
  \end{letlist}
\end{theorem}

This theorem includes as special cases  a number of known results for 
homogeneous and inhomogenous percolation on the square,
triangular, and hexagonal lattices beginning with Kesten's theorem that $\pc=\frac12$
for the square lattice, see \cite{Kesten80,Kesten_book,ZSWS}.

We turn now to the universality of critical exponents. Recall the exponents $\rho$, $\eta$, and $\de$ 
of Table \ref{Tab-ce}. 
The exponent $\rho_{2j}$ is the so-called $2j$ alternating-arm critical exponent, see \cite{G-rev,GM3}. 
An exponent is said to be $\sG$-\emph{invariant} if its value
is constant across the family $\sG$.

\begin{theorem}[Universality \cite{GM3}]\label{thm:2}
Let $\sG$ be the class of isoradial graphs with the bounded-angles property and the square-grid property.
\begin{letlist}
\item
Let $\pi \in\{\rho\}\cup\{\rho_{2j}: j \ge 1\}$.
If $\pi$ exists for some $G \in  \sG$, then it is $\sG$-invariant.
\item
If either $\rho$ or $\eta$ exists for some $G \in \sG$, then $\rho$, $\eta$, $\de$ are
$\sG$-invariant and satisfy the scaling relations $\eta\rho=2$ and $2\rho = \de+1$.
\end{letlist}
\end{theorem}

The theorem establishes universality for bond percolation on isoradial graphs,
but restricted to the exponents $\rho$, $\eta$, $\de$ \emph{at the critical point}. The method of proof
does not seem to extend to the near-critical exponents $\beta$, $\gamma$, etc (see Problem E of Section
\ref{sec:open-p}). 

It is in fact `known' that, for reasonable two-dimensional lattices,
\begin{equation}\label{exact}
\rho=\tfrac{48}5, \quad \eta=\tfrac5{24}, \quad\de=\tfrac{91}5,
\end{equation}
although these values (and more), 
long predicted in the physics literature, have been proved rigorously
only for (essentially) site percolation on the triangular lattice. See Lawler, Schramm, Werner \cite{LSW01}
and Smirnov and Werner \cite{Smirnov-Werner}. Note that site percolation
on the triangular lattice does not lie within the ambit of
Theorems \ref{thm:1} and \ref{thm:2}.   

To summarize, there is currently no
known proof of the existence of critical exponents for any graph belonging to $\sG$. However, if
certain exponents exist for \emph{any} such graph, then they exist for all $G$ and are $\sG$-invariant.
If one could establish a result such as in \eqref{exact} for any such graph, 
then this result would be valid across the entire family $\sG$.

The main ideas of the proofs of Theorems \ref{thm:1} and \ref{thm:2} are as follows.
The first element is the so-called \bxp.
Loosely speaking, this is the property that the probability of
an open crossing of a box with given aspect-ratio is bounded away from 0, uniformly in the
position, orientation, and size of the box.
The \bxp\ was proved by Russo \cite{Russo} and Seymour/Welsh \cite{Seymour-Welsh}  
for homogeneous percolation on
the square lattice, using its properties of symmetry and self-duality.
It may be shown using classical methods that the
\bxp\ is a certificate of a critical or supercritical percolation
model. It may be deduced that, if both the primal and dual models have the \bxp, then
they are both critical. 

The \stt\ of the next section provides a method for transforming one isoradial graph
into another. The key step in the proofs is to show that this transformation preserves the \bxp.
It follows that any isoradial graph that can be obtained by a sequence of
transformations from the square lattice has the \bxp,
and is therefore critical. It is proved in \cite{GM3} that this includes any isoradial
graph with both the BAP and SGP.
 
\subsection{Star--triangle transformation}\label{sec:stt}

\begin{figure}
 \begin{center}
 \includegraphics[width=0.8\textwidth]{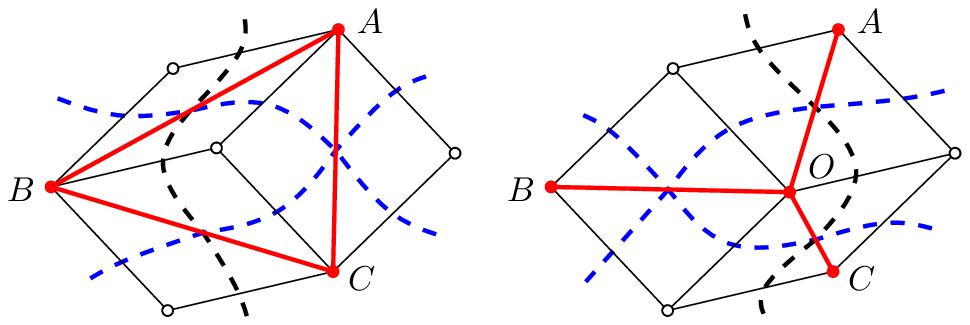}\newline
 \includegraphics[width=0.7\textwidth]{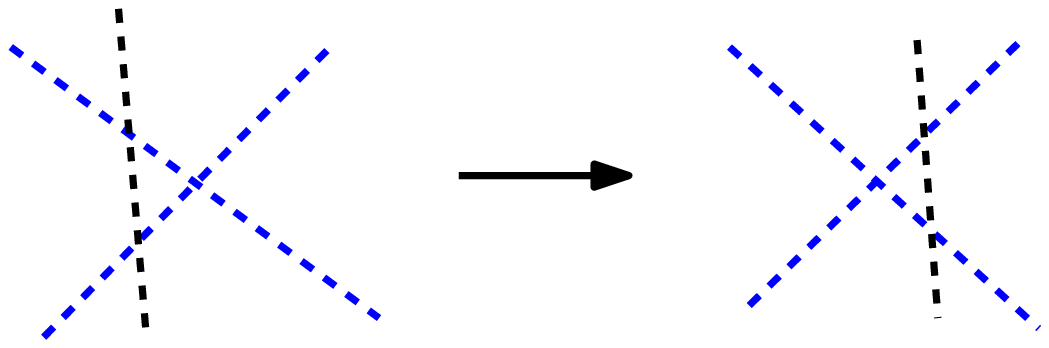}
\end{center}
 \caption{There are two ways of tiling the hexagon in the upper figure,
and switching between these amounts to a \stt\ for the isoradial graph.
The effect on the track system is illustrated in the lower figure.}
 \label{fig:stt-on-iso}
\end{figure}

The central fact that permits proofs of criticality and universality
is that the \stt\ has a geometric representation that acts locally on
rhombic tilings. Consider three rhombi assembled to form a hexagon as in 
the upper left of Figure \ref{fig:stt-on-iso}. The interior of the hexagon may be 
tiled by (three) rhombi in either of two ways, the other such tiling being drawn at the upper right of
the figure. The switch from the first to the second has two effects:
(i) the track system is altered as indicated there, with one track being moved over the
intersection of the other two, and (ii) the triangle in the isoradial graph of the upper left
is transformed into a star.  These observations are graph-theoretic rather than model-specific.
We next parametrize the system in such a way that the parameters mutate
in the canonical way under the above transformation. That is, for a given probabilistic model,  we seek
a parametrization under which the geometrical switch induces the appropriate
parametric change.

\begin{figure}
 \centering
	\psfrag{A}{$A$}
	\psfrag{B}{$B$}
	\psfrag{C}{$C$}
	\psfrag{O}{$O$}
	\psfrag{p0}{$p_0$}
	\psfrag{p1}{$p_1$}
	\psfrag{p2}{$p_2$}
	\psfrag{p0p}{$1-p_0$}
	\psfrag{p1p}{$1-p_1$}
	\psfrag{p2p}{$1-p_2$}
    \includegraphics[width=0.6\textwidth]{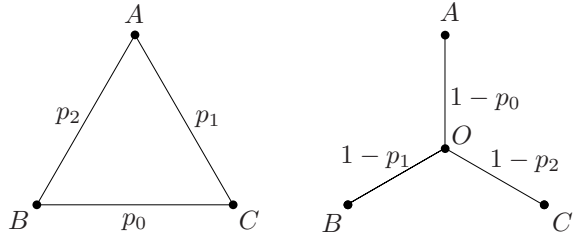}
  \caption{The star--triangle transformation for bond percolation.}
  \label{fig:star_triangle_transformation}
\end{figure}

Here is the \stt\ for percolation. 
Consider the triangle $T=(V,E)$
and the star $S=(V',E')$ as drawn in Figure~\ref{fig:star_triangle_transformation}. Let $\bp=(p_0,p_1,p_2)\in[0,1)^3$, and suppose the edges in the figure are declared
open with the stated probabilities.
The two ensuing configurations induce two connectivity relations on the set $\{A,B,C\}$ within
$S$ and $T$, \resp. It turns out that these two connectivity relations are equi-distributed so long
as $\kappa(\bp)=0$, where
\begin{equation}\label{kappadef}
\kappa(\bp) = p_0+p_1+p_2 - p_1p_2p_3 -1.
\end{equation}

The \stt\ is used as follows. Suppose, in a graph $G$, one finds a triangle whose
edge-probabilities
satisfy \eqref{kappadef}. Then this triangle may be replaced by a star having the
complementary probabilities of Figure \ref{fig:star_triangle_transformation} 
without altering the probabilities of any
long-range connections in $G$. Similarly, stars may be transformed into triangles. One complicating
feature of the transformation is the creation of a new vertex when passing from
a triangle to a star (and its destruction when passing in the reverse direction). 

The \stt\  was discovered first in the context of electrical networks
by Kennelly \cite{Ken} in 1899, and
it was adapted in 1944 by Onsager \cite{OnsI} to the Ising model in
conjunction with Kramers--Wannier duality.
It is a key element in the work
of Baxter \cite{Bax78,Baxter_book} on exactly solvable models in statistical mechanics, where it has
become known as the \emph{Yang--Baxter equation} (see \cite{Perk-AY}
for a history of its importance in physics).
Sykes and Essam \cite{Sykes_Essam} used the \stt\
to predict the critical surfaces of certain inhomogeneous (but periodic) bond percolation 
processes on the
triangular and hexagonal lattices, and furthermore the \stt\  is a tool in the study of
the \rc\ model \cite[Sect.\ 6.6]{G-rcm}, and the dimer model \cite{BdT11}.

Let us now explore the operation of the \stt\ in the context
of the rhombic switch of Figure \ref{fig:stt-on-iso}. Let
$G$ be an isoradial graph containing the upper left hexagon of the figure, and let $G'$ be the new graph after
the rhombic switch.
The definition \eqref{pdef}
of the edge-probabilities has been chosen in such a way that the
values on the triangle satisfy \eqref{kappadef} and those on the star are as given in
Figure \ref{fig:star_triangle_transformation}. It follows that the connection
probabilities on $G$ and $G'$ are equal. 
Graphs which have been thus parametrized but not embedded isoradially were
called \emph{$Z$-invariant} by Baxter \cite{Bax78}.
See \cite{LT} for a recent account of the application of the above rhombic switch to Glauber dynamics of 
lozenge tilings of the triangular lattice.

One may couple the probability spaces on $G$ and $G'$ in such a way that the 
\stt\ preserves \emph{open connections}, rather than just their probabilities.
Suppose that, in a given configuration, there exists an open path in $G$ between vertex-sets $A$ and $B$.
On applying a sequence of \stt s, we obtain an open path in $G'$ from
the image of $A$ to the image of $B$. Thus, \stt s transport open paths to open paths, and it is thus
that the \bxp\ is transported from $G$ to $G'$.

In practice, infinitely many \stt s are required to achieve the necessary
transitions between graphs. The difficulties of the proofs of Theorems \ref{thm:1}--\ref{thm:2}
are centred on the need to establish 
sufficient control on the drifts of paths and their endvertices under these transformations.

\section{Open problems for percolation}\label{sec:open-p}

We discuss associated open problems in this section.

\begin{problist}

\item[A.] 
\emph{Existence and equality of critical exponents.}
It is proved in Theorem \ref{thm:2} that, if the three exponents $\rho$, $\eta$, $\de$ exist
for some member of the family $\sG$, then they exist for all members of
the family, and are constant across the family.
Essentially the only model for which existence has been proved is the site model on the
triangular lattice, but this does not belong to $\sG$. A proof of existence of exponents for the bond model on
the square lattice would imply their existence for the isoradial graphs studied
here. Similarly, if one can show any
exact value for the latter bond model, then this value holds across $\sG$.

\item[B.]
\emph{Cardy's formula.}
Smirnov's proof \cite{Smirnov} of Cardy's formula  has resisted
extension to models beyond  site percolation on the triangular lattice.
It seems likely that Cardy's formula is valid for canonical percolation on
any reasonable isoradial graph. There is a strong sense in which the existence of interfaces  is
preserved under the \stt s of the proofs. On the other hand,  there is currently only limited control of
the geometrical perturbations of interfaces, and in addition Cardy's formula is as 
yet unproven for 
\emph{all} isoradial bond percolation models.

\item[C.]
\emph{The bounded-angles property.}
It is normal in working with probability and isoradial graphs
to assume the BAP, see for example \cite{Chelkak-Smirnov3}. In the language of
finite element methods, \cite{JKK}, the BAP is an example of the \v Zen\'\i \v sek--Zl\'amal condition.

The BAP is a type of uniform non-flatness assumption. It implies an equivalence
of metrics, and enables a uniform boundedness of certain probabilities.
It may, however,  not be \emph{necessary} for the \bxp, and hence for the main results above.

As a test case, consider the situation in which all rhombi have angles exactly $\eps$ and $\pi-\eps$. In the limit
as $\eps\downarrow 0$, we obtain\footnote{Joint work with Omer Angel.} 
the critical space--time percolation process
on $\ZZ \times \RR$, see Figure \ref{fig:cts} and, for example, \cite{G-st}. 
Let $B_n(\alpha)$ be an $n \times n$ square
of $\RR^2$ inclined at angle $\a$, and let $C_n(\alpha)$ be the event that the square is traversed
by an open path between two given opposite faces. It is elementary using duality that
$$
\PP\bigl(C_n(\tfrac14\pi)\bigr) \to \tfrac12 \qq\text{as } n \to \oo.
$$
Numerical simulations (of A.\ Holroyd) suggest that the same limit holds when $\a=0$.
A proof of this would suggest that the limit does not depend on $\a$, and this in turn 
would support the 
possibility that the critical space--time percolation process satisfies Cardy's formula.

\begin{figure}
  \begin{center}
	\psfrag{L}{$L$}
    \includegraphics[width=0.45\textwidth]{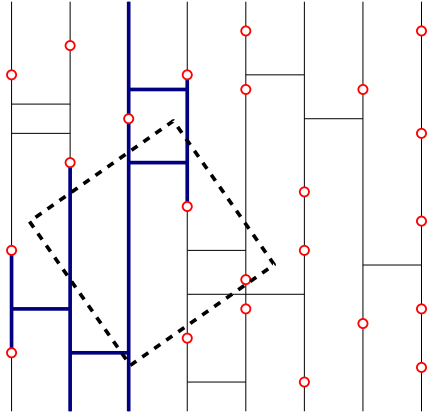}
  \end{center}
  \caption{\emph{Space--time percolation.} Each line is cut at rate $1$,
and nearest neighbours are joined at rate $1$. One of the open clusters is highlighted. We
ask for the probability that the box is traversed by an open path from its lower left side to its upper right side.}
  \label{fig:cts}
\end{figure}

\item[D.]
\emph{The square-grid property.} 
The SGP is a useful tool in the proof of Theorem \ref{thm:2},
but it may not be necessary. In \cite{GM3} is presented an isoradial graph without the SGP,
and this example may be handled using an additional ad hoc argument.

\item[E.]
\emph{Near-critical exponents.}
Theorem \ref{thm:2} establishes the universality of exponents \emph{at criticality}.
The method of proof does not appear to be extendable to the \emph{near-critical
exponents}, and it is an open problem to prove these to be universal for isoradial graphs.
Kesten showed in \cite{Kesten87} (see also \cite{Nolin}) that certain properties of a critical percolation
process imply properties of the near-critical process, so long as the underlying graph
has a sufficiently rich automorphism group. In particular, 
for such graphs, knowledge of certain critical
exponents at criticality implies knowledge of exponents away from criticality.
Only certain special isoradial graphs have sufficient homogeneity for such arguments
to hold without new ideas of substance, and it is an open problem to weaken 
these assumptions of homogeneity. See the discussion around \cite[Thm 1.2]{GM2}.

\item[F.]
\emph{Random-cluster models.}
How far may the proofs be extended to other models? It may seem
at first sight that only a \stt\ is required, but, as usual in such situations,
boundary conditions play a significant role for dependent models such
as the \rc\ model. The control of boundary conditions presents a new
difficulty, so far unexplained.  We return to this issue in Section \ref{sec:open-rcm}.

\end{problist}

\section{Random-cluster model}\label{sec:rcm}

\subsection{Background}

The \rc\ model was introduced by Fortuin and Kasteleyn around 1970 as a unification
of processes satisfying versions of the series and parallel laws.
In its base form, the \rc\ model has two parameters, an edge-parameter
$p$ and a cluster-weighting factor $q$. 

Let $G=(V,E)$ be a finite graph, with associated configuration space   $\Om=\{ 0,1\}^E$. 
For $\om\in\Om$ and $e \in E$, the edge $e$ is designated \emph{open}
if $\om_e=1$. Let $k(\om)$ be the number of open clusters
of a configuration $\om$.
The \emph{\rc\ measure\/} on $\Om$,
with parameters $p\in[0,1]$, $q\in(0,\oo)$, is the probability measure
satisfying
\begin{equation}\label{rcmeas}
\fpq (\om )\propto  q^{k(\om )}\PP_p(\om) ,\qq\om\in\Om,
\end{equation}
where $\PP_p$ is the percolation product-measure with density $p$.
In a more general setting, each edge $e \in E$ has an associated parameter $p_e$.

Bond percolation is retrieved by setting $q=1$,
and electrical networks
arise via the limit $p,q\to 0$ in such a way that $q/p\to 0$.
The relationship to Ising/Potts models is more complicated
and involves a transformation
of measures. In brief, two-point connection
probabilities for the \rc\ measure with $q\in\{2,3,\dots\}$
correspond to two-point correlations for ferromagnetic $q$-state Ising/Potts models,
and this allows a geometrical interpretation
of the latter's correlation structure.
A fuller account of the \rc\ model and its history and associations may be found in \cite{G-rcm,Werner_SMF},
to which the reader is referred for the basic properties of the model.

The special cases of percolation and the Ising model are very much better understood
than is the general \rc\ model. We restrict ourselves to two-dimensional
systems in this review, and we concentrate on the question of the identification of critical surfaces
for certain isoradial graphs. 

Two pieces of significant recent progress are reported here. Firstly, Beffara and Duminil-Copin
\cite{Beffara_Duminil} have developed the classical approach of percolation in order to
identify the critical point of the square lattice, thereby solving a longstanding conjecture.
Secondly, together with Smirnov \cite{BDS2}, they have made use of the so-called parafermionic 
observable of \cite{Smi10} in a study of the critical surfaces of periodic isoradial graphs with $q \ge 4$.

\subsection{Formalities}
The \rc\ measure may not be defined directly on an \emph{infinite} graph $G$.
There are two possible ways to proceed in the setting of an infinite graph, 
namely via either boundary conditions or the DLR condition.
The former approach works as follows. Let $(G_n: n \ge 1)$ be an increasing
family of finite subgraphs of $G$ that exhaust $G$ in the limit $n \to \oo$, and let
$\pd G_n$ be the \emph{boundary} of $G_n$, that is,
$\pd G_n$ is  the set of vertices of $G_n$ that are adjacent to a vertex of $G$ not in $G_n$.
A \emph{boundary condition} is an equivalence relation  $b_n$ on $\pd G_n$;
any two vertices $u,v \in \pd G_n$ that are equivalent under $b_n$ are taken to be part of the same
cluster. The extremal boundary conditions are: the \emph{free} boundary
condition, denoted $b_n = 0$, for which each vertex is in a separate
equivalence class; and the \emph{wired} boundary condition, denoted 
$b_n = 1$, with a unique equivalence class.
We now consider the
set of weak limits as $n \to \oo$ of the \rc\ measures on $G_n$ with boundary conditions $b_n$.

Assume henceforth that  $q \ge 1$. Then the \rc\ measures have properties of 
positive association and stochastic ordering, and one may deduce that the
free and wired boundary conditions $b_n= 0$ and $b_n = 1$ are extremal in the following sense:
(i) there is a unique weak limit of the free measures 
(\resp, the wired measures),
and (ii)  any other weak limit lies, in the sense of stochastic ordering, between
these two limits. We write $\fpq^0$ and $\fpq^1$ for the free and wired weak limits.
It is an important question to determine when $\fpq^0=\fpq^1$, and the
answer so far is incomplete even when $G$ has a periodic structure,
see \cite[Sect.\ 5.3]{G-rcm}.

The \emph{percolation probabilities} are defined by
\begin{equation}\label{13.8}
\theta^b(p,q)=\fpqb (0\lra\infty ),\qq b=0,1,
\end{equation}
and the \emph{critical values\/} by
\begin{equation}\label{rccritprob}
\pcb (q)=\sup\{ p:\theta^b(p,q)=0\} ,\qq b=0,1.
\end{equation}

Suppose that $G$ is embedded in $\RR^d$ in a natural manner. 
When $G$ is \emph{periodic} (that is, its embedding is invariant under a $\ZZ^d$ action),
there is a general argument using convexity of pressure 
(see \cite{G93}) that implies that $\pc^0(q)=\pc^1(q)$, and in this case we write
$\pc(q)$ for the common value. 

One of the principal problems is to determine for which $q$ the percolation probability
$\theta^1(p,q)$ is discontinuous at the critical value $\pc$. This amounts to asking when
$\theta^1(\pc,q)>0$;  the phase transition is said to be of \emph{first order} whenever the last 
inequality holds. 
The phase transition is known to be of first order for sufficiently large $q$, and is
believed to be so if and only if $q>Q(d)$ for some $Q(d)$ depending on the dimension $d$. 
Furthermore, it is expected that
$$
Q(d) = \begin{cases} 4 &\text{if } d=2,\\
2&\text{if } d \ge 4.
\end{cases}
$$
We restrict our attention henceforth to the case $d=2$, for which 
it is believed that the value $q=4$ separates the first and second order transitions. 
Recall Conjecture \ref{thetapc0} and note the recent proof that
$Q(2) \ge 4$, for which the reader is referred to
 \cite{DC-book} and the references therein.

\subsection{Critical point on the square lattice}\label{sec:critpt}

The square lattice $\ZZ^2$ is one of the main playgrounds of physicists and probabilists. Although the
critical points of percolation, the Ising model and some Potts models on $\ZZ^2$
are long proved, the general answer for \rc\ models (and hence all Potts models) has been proved
only recently.

\begin{theorem}[Criticality \cite{Beffara_Duminil}]\label{thm@pc}
The \rc\ model on the square lattice with cluster-weighting factor $q \ge 1$ has
critical value
$$
\pc(q) = \frac {\sqrt q}{1+\sqrt q}.
$$
\end{theorem}

This exact value has been `known' for a long time.
When $q=1$, the statement $\pc(1)=\frac12$ is the Harris--Kesten theorem for bond percolation.
When $q=2$, it amounts to the well known
calculation of the critical temperature of the
Ising model. For large $q$, the result (and more) was proved  in \cite{LMMRS, LMR}
($q> 25.72$ suffices, see \cite[Sect. 6.4]{G-rcm}).
There is a `physics proof' in  \cite{HKW} for $q \ge 4$.

\begin{figure}
  \begin{center}
	\psfrag{L}{$L$}
    \includegraphics[width=0.6\textwidth]{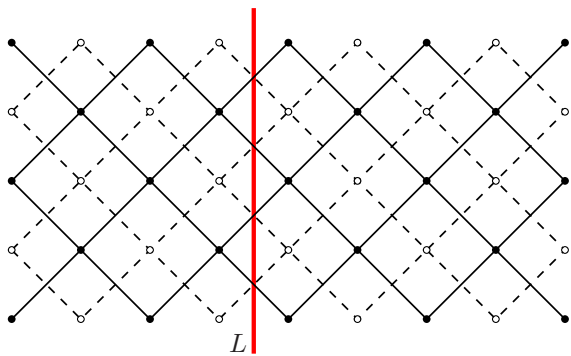}
  \end{center}
  \caption{The square lattice and its dual, rotated through $\pi/4$. Under reflection in the line $L$,
the primal is mapped to the dual.}
  \label{fig:z2rot}
\end{figure}

The main contribution of \cite{Beffara_Duminil} is a proof of a \bxp\
using a clever extension
of the `RSW' arguments of Russo and Seymour--Welsh in the context of
the symmetry illustrated in Figure \ref{fig:z2rot}, combined with careful control of boundary
conditions. An alternative approach is developed in \cite{DM14}.

\subsection{Isoradiality and the \stt}
\begin{figure}[t]
  \begin{center}
	\psfrag{A}{$A$}
	\psfrag{B}{$B$}
	\psfrag{C}{$C$}
	\psfrag{O}{$O$}
	\psfrag{y0}{$y_0$}
	\psfrag{y1}{$y_1$}
	\psfrag{y2}{$y_2$}
	\psfrag{y0p}{$y_0'$}
	\psfrag{y1p}{$y_1'$}
	\psfrag{y2p}{$y_2'$}
    \includegraphics[width=0.6\textwidth]{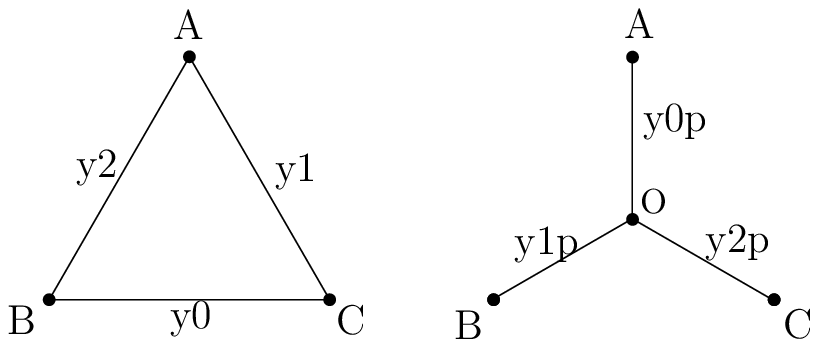}
  \end{center}
  \caption{The \stt\ for the \rc\ model.}
  \label{fig:stt-rcm}
\end{figure}

The \stt\ for the \rc\ model is similar to that of percolation, and is illustrated in Figure \ref{fig:stt-rcm}.
The three edges of the triangle have parameters $p_0$, $p_1$, $p_2$,
and we set $\mathbf y =(y_0,y_1,y_2)$ where
$$
y_i= \frac{p_i}{1-p_i}.
$$
The corresponding edges of the star have parameters $y_i'$ where $y_iy_i'=q$.
Finally, we require that the $y_i$ satisfy $\psi(\mathbf y)=0$ where
\begin{equation}\label{kappadef2}
\psi(\mathbf y) = y_1y_2y_3 +y_1y_2 + y_2y_3 + y_3y_1-q.
\end{equation}
Further details of the \stt\ for the \rc\ model may be found in \cite[Sect.\ 6.6]{G-rcm}.

We now follow the discussion of Section \ref{sec:stt} of the relationship between the \stt\ and
the rhombus-switch of Figure \ref{fig:stt-on-iso}. In so doing, we arrive 
(roughly as in \cite[p.\ 282]{Ken02}) at the `right' parametrization
for an isoradial graph $G$, namely with \eqref{pdef}  replaced by
\begin{equation}
\begin{alignedat}{2}\label{pdef3}
\text{if $1\le q < 4$:}\q y_e &= 
\sqrt q\,\dfrac{\sin(\tfrac12\sigma(\pi-\theta_e))}{\sin\bigl(\tfrac12\sigma\theta_e\bigr)},&\q
\cos(\tfrac12\sigma\pi) &= \tfrac12\sqrt q,\\
\text{if $q > 4$:}\q y_e&= 
\sqrt q\,\dfrac{\sinh(\tfrac12\sigma(\pi-\theta_e))}{\sinh\bigl(\tfrac12\sigma\theta_e\bigr)},\q&
\cosh(\tfrac12\sigma\pi) &= \tfrac12\sqrt q,
\end{alignedat}
\end{equation}
where $\theta_e$ is given in Figure \ref{fig:isomap}.
The intermediate case $q=4$ is the common limit of the two expressions as $q \to4$, namely
$$
y_e = 2\,\frac{\pi-\theta_e}{\theta_e}.
$$

Write $\fgq^b$ for the corresponding \rc\ measure on an isoradial graph $G$ with boundary condition $b=0,1$.
We refer to $\fgq^0$ as the `canonical \rc\ measure' on $G$.

\subsection{Criticality via the parafermion}\label{sec:paraf}
Theorem \ref{thm@pc} is proved in \cite{Beffara_Duminil} by classical methods, 
and it holds for all $q\ge 1$.
The proof is sensitive to the assumed symmetries of the lattice, and does not currently extend
even to the inhomogeneous \rc\ model on $\ZZ^2$ in which the vertical and horizontal edges
have different parameter values. In contrast, the parafermionic observable introduced
by Smirnov \cite{Smi10} has been exploited by Beffara, Duminil-Copin, and Smirnov \cite{BDS2}
to study the critical point of fairly general isoradial graphs subject to the condition $q \ge 4$.

Let $G=(V,E)$ be an isoradial graph.
For $\b\in(0,\oo)$, let $y_e(\b)=\b y_e$ where
$y_e$ is given in \eqref{pdef3}. Let 
$$
p_e(\b)=\frac{y_e(\b)}{1+y_e(\b)}
$$ 
accordingly, and write
$\fgbq^b$ for the corresponding \rc\ measure on $G$ with boundary condition $b$. 
The following result of \cite{BDS2} is proved by a consideration of the parafermionic observable.

\begin{theorem}[\cite{BDS2}]\label{thm:bds0}
Let $q \ge 4$, and let $G$ be an isoradial graph satisfying the BAP.
For $\b<1$, there exists $a>0$ such that
$$
\fgbq^0(u \lra v) \le e^{-a|u-v|}, \qq u,v \in V.
$$
\end{theorem}

One deduces from Theorem \ref{thm:bds0} using duality that
\begin{letlist}
\item for $\b<1$, $\fgbq^{0}$-a.s., there is  no infinite open cluster, and
\item for $\b>1$, $\fgbq^{1}$-a.s., there exists a unique infinite open cluster.
\end{letlist}
This is only a partial verification of the criticality of the canonical measure,
since parts (a) and (b) deal with potentially different measures, 
namely the free and wired limit measures,
\resp. Further progress may be made for \emph{periodic} graphs, as follows. 
Subject to the assumption of periodicity,
it may be proved as in \cite{G93} that $\fgbq^0=\fgbq^1$ for almost every
$\b$, and hence that part (b) holds with $\fgbq^{1}$
replaced by $\fgbq^{0}$. Therefore, for periodic embeddings,
the canonical measure $\fgq^0=\phi^0_{G,q,1}$ is critical.
 
Here is an application of the above remarks to the (periodic) inhomogeneous square lattice.

\begin{corollary}[\cite{BDS2}]\label{cor:sq}
Let $q \ge 4$, and consider the \rc\ model on $\ZZ^2$ with the variation that
horizontal edges have parameter $p_1$ and
vertical edges parameter $p_2$.
The critical surface is given by $y_1y_2=q$ where
$y_i=p_i/(1-p_i)$.
\end{corollary}

We close with the observation that a great deal more is known in
the special case when $q=2$. The $q=2$ \rc\ model corresponds to the Ising model, for which the
special arithmetic of the equation $1+1=2$ permits a number of
techniques which are not available in greater generality. In particular, the Ising model
and the $q=2$ \rc\ model
on an isoradial graph lend themselves to a fairly complete theory using the parafermionic observable.
The interested reader is directed to the work of Smirnov \cite{Smi07,Smi10} and Chelkak--Smirnov 
\cite{Chelkak-Smirnov2}.

\section{Open problems for the \rc\ model}\label{sec:open-rcm}

\begin{problist}

\item[A.]
\emph{Inhomogeneous models.}
Extend Corollary \ref{cor:sq} to cover the case $1\le q < 4$. 

\item[B.]
\emph{Periodicity.}
Remove the assumption of periodicity in the proof of criticality of the canonical \rc\ measure
on isoradial graphs.
It would suffice to prove that $\fgbq^0=\fgbq^1$ for almost every $\b$, without the assumption
of periodicity. More generally, it would be useful to have a proof of the uniqueness of Gibbs states
for \emph{aperiodic} interacting systems, along the lines of that of Lebowitz and Martin-L{\"o}f \cite{LebML}
for a \emph{periodic} Ising model.

\item[C.] 
\emph{Bounded-angles property.}
Remove the assumption of the bounded-angles property in Theorem \ref{thm@pc}.

\item[D.]
\emph{Criticality and universality for general $q$.}
Adapt the arguments of \cite{GM3} (or otherwise) to prove criticality and universality for
the canonical \rc\ measure on isoradial graphs either for general $q \ge 1$ or subject 
to the restriction $q \ge 4$.

\end{problist}

\section*{Acknowledgements}
The author is grateful to Ioan Manolescu for many discussions concerning 
percolation on isoradial graphs, and to Omer Angel and Alexander Holroyd for
discussions about the space--time percolation process of Figure \ref{fig:cts}. 
Hugo Duminil-Copin and Ioan Manolescu kindly commented on a draft of this paper.
This work was supported in part by the EPSRC under grant EP/103372X/1.

\bibliography{icm}

\providecommand{\bysame}{\leavevmode\hbox to3em{\hrulefill}\thinspace}
\providecommand{\MR}{\relax\ifhmode\unskip\space\fi MR }
\providecommand{\MRhref}[2]{%
  \href{http://www.ams.org/mathscinet-getitem?mr=#1}{#2}
}
\providecommand{\href}[2]{#2}
\begin{thebibliography}{10}

\bibitem{BDGS}
R.~Bauerschmidt, H.~Duminil-Copin, J.~Goodman, and G.~Slade, \emph{Lectures on
  self-avoiding-walks}, Probability and Statistical Physics in Two and More
  Dimensions (D.~Ellwood, C.~M. Newman, V.~Sidoravicius, and W.~Werner, eds.),
  Clay Mathematics Institute Proceedings, vol.~15, CMI/AMS publication, 2012,
  pp.~395--476.

\bibitem{Bax78}
R.~J. Baxter, \emph{Solvable eight-vertex model on an arbitrary planar
  lattice}, Philos. Trans. Roy. Soc. London Ser. A \textbf{289} (1978),
  315--346.

\bibitem{Baxter_book}
\bysame, \emph{{Exactly Solved Models in Statistical Mechanics}}, Academic
  Press, London, 1982.

\bibitem{Beffara_Duminil}
V.~Beffara and H.~{Duminil-Copin}, \emph{The self-dual point of the
  two-dimensional random-cluster model is critical for $q \geq 1$}, Probab. Th.
  Rel. Fields \textbf{153} (2012), 511--542.

\bibitem{BDS2}
V.~Beffara, H.~{Duminil-Copin}, and S.~Smirnov, \emph{On the critical
  parameters of the $q\ge 4$ random-cluster model on isoradial graphs},
  (2013), preprint.

\bibitem{BolRio}
B.~Bollob\'as and O.~Riordan, \emph{Percolation}, Cambridge University Press,
  Cambridge, 2006.

\bibitem{BdT11}
C.~Boutillier and B.~de Tili\`ere, \emph{The critical {$Z$}-invariant {I}sing
  model via dimers: Locality property}, Commun. Math. Phys. \textbf{301}
  (2011), 473--516.

\bibitem{BdT12}
\bysame, \emph{Statistical mechanics on isoradial graphs}, {Probability in
  Complex Physical Systems} (J.-D. Deuschel, B.~Gentz, W.~K\"onig, M.~von
  Renesse, M.~Scheutzow, and U.~Schmock, eds.), Springer Proceedings in
  Mathematics, vol.~11, 2012, pp.~491--512.

\bibitem{JKK}
J.~Brandts, S.~Korotov, and M.~K{\v r}\'\i{\v z}ek, \emph{Generalization of the
  {Z}l\' amal condition for simplicial finite elements in {$\RR^d$}}, Applic.
  Math. \textbf{56} (2011), 417--424.

\bibitem{BrH}
S.~R. Broadbent and J.~M. Hammersley, \emph{Percolation processes {I. C}rystals
  and mazes}, Proc. Camb. Phil. Soc. \textbf{53} (1957), 629--641.

\bibitem{deB1}
N.~G.~de Bruijn, \emph{Algebraic theory of {P}enrose's non-periodic tilings of
  the plane. {I}}, Indagat. Math. (Proc.) \textbf{84} (1981), 39--52.

\bibitem{deB2}
\bysame, \emph{Algebraic theory of {P}enrose's non-periodic tilings of the
  plane. {II}}, Indagat. Math. (Proc.) \textbf{84} (1981), 53--66.

\bibitem{deB3}
\bysame, \emph{Dualization of multigrids}, J. Phys. Colloq. \textbf{47} (1986),
  85--94.

\bibitem{Cardy}
J.~Cardy, \emph{Critical percolation in finite geometries}, J. Phys. A: Math.
  Gen. \textbf{25} (1992), L201--L206.

\bibitem{Chelkak-Smirnov3}
D.~Chelkak and S.~Smirnov, \emph{Discrete complex analysis on isoradial
  graphs}, Adv. Math. \textbf{228} (2011), 1590--1630.

\bibitem{Chelkak-Smirnov2}
\bysame, \emph{{Universality in the 2D Ising model and conformal invariance of
  fermionic observables}}, Invent. Math. \textbf{189} (2012), 515--580.

\bibitem{Duff}
R.~J. Duffin, \emph{Potential theory on a rhombic lattice}, J. Combin. Th.
  \textbf{5} (1968), 258--272.

\bibitem{DC-book}
H.~Duminil-Copin, \emph{Parafermionic observables and their applications to
  planar statistical physics models}, Ensaios Matem\'aticos \textbf{25} (2013),
  1--371.

\bibitem{DM14}
H.~Duminil-Copin and I.~Manolescu, \emph{The phase transitions of the planar
  random-cluster model and {P}otts model with $q \ge 1$ is sharp},  (2014), in
  preparation.

\bibitem{Fitz}
R.~J. Fitzner, \emph{Non-backtracking lace expansion}, Ph.D. thesis, Technische
  Universiteit Eindhoven, 2013.

\bibitem{G93}
G.~R. Grimmett, \emph{The stochastic random-cluster process and the uniqueness
  of random-cluster measures}, Ann. Probab. \textbf{23} (1995), 1461--1510.

\bibitem{G99}
\bysame, \emph{Percolation}, 2nd ed., Springer, Berlin, 1999.

\bibitem{G-rcm}
\bysame, \emph{{The Random-Cluster Model}}, Springer, Berlin, 2006.

\bibitem{G-st}
\bysame, \emph{Space--time percolation}, In and Out of Equilibrium 2
  (V.~Sidoravicius and M.~E. Vares, eds.), Progress in Probability, vol.~60,
  Birkh\"auser, Boston, 2008, pp.~305--320.

\bibitem{Grimmett_Graphs}
\bysame, \emph{{Probability on Graphs}}, Cambridge University Press, Cambridge,
  2010, \url{http://www.statslab.cam.ac.uk/~grg/books/pgs.html}.

\bibitem{G-rev}
\bysame, \emph{Three theorems in discrete random geometry}, Probab. Surveys
  \textbf{8} (2011), 403--441.

\bibitem{GM1}
G.~R. Grimmett and I.~Manolescu, \emph{Inhomogeneous bond percolation on the
  square, triangular, and hexagonal lattices}, Ann. Probab. \textbf{41} (2013),
  2990--3025.

\bibitem{GM2}
\bysame, \emph{Universality for bond percolation in two dimensions}, Ann.
  Probab. \textbf{41} (2013), 3261--3283.

\bibitem{GM3}
\bysame, \emph{Bond percolation on isoradial graphs: criticality and
  universality}, Probab. Th. Rel. Fields (2014),
  \url{http://arxiv.org/abs/1204.0505}.

\bibitem{H57a}
J.~M. Hammersley, \emph{Percolation processes. {L}ower bounds for the critical
  probability}, Ann. Math. Statist. \textbf{28} (1957), 790--795.

\bibitem{H59}
\bysame, \emph{Bornes sup\'erieures de la probabilit\'e critique dans un
  processus de filtration}, {Le Calcul des Probabilit\'es et ses Applications},
  CNRS, Paris, 1959, pp.~17--37.

\bibitem{HS}
T.~Hara and G.~Slade, \emph{Mean-field critical behaviour for percolation in
  high dimensions}, Commun. Math. Phys. \textbf{128} (1990), 333--391.

\bibitem{HS94}
\bysame, \emph{Mean-field behaviour and the lace expansion}, {Probability and
  Phase Transition} (G.~R. Grimmett, ed.), Kluwer, 1994, pp.~87--122.

\bibitem{HKW}
D.~Hintermann, H.~Kunz, and F.~Y. Wu, \emph{Exact results for the {P}otts model
  in two dimensions}, J. Statist. Phys. \textbf{19} (1978), 623--632.

\bibitem{Ken}
A.~E. Kennelly, \emph{The equivalence of triangles and three-pointed stars in
  conducting networks}, Electrical World and Engineer \textbf{34} (1899),
  413--414.

\bibitem{Ken02}
R.~Kenyon, \emph{An introduction to the dimer model}, School and Conference on
  Probability Theory (G.~F. Lawler, ed.), Lecture Notes Series, vol.~17, ICTP,
  Trieste, 2004, \url{http://publications.ictp.it/lns/vol17/vol17toc.html},
  pp.~268--304.

\bibitem{KenS}
R.~Kenyon and J.-M. Schlenker, \emph{Rhombic embeddings of planar quad-graphs},
  Trans. Amer. Math. Soc. \textbf{357} (2005), 3443--3458.

\bibitem{Kesten80}
H.~Kesten, \emph{The critical probability of bond percolation on the square
  lattice equals $1/2$}, Commun. Math. Phys. \textbf{74} (1980), 44--59.

\bibitem{Kesten_book}
\bysame, \emph{{Percolation Theory for Mathematicians}}, Birkh\"auser, Boston,
  1982.

\bibitem{Kesten87}
\bysame, \emph{Scaling relations for {2D}-percolation}, Commun. Math. Phys.
  \textbf{109} (1987), 109--156.

\bibitem{KN}
G.~Kozma and A.~Nachmias, \emph{Arm exponents in high dimensional percolation},
  J. Amer. Math. Soc. \textbf{24} (2011), 375--409.

\bibitem{LMMRS}
L.~Laanait, A.~Messager, S.~Miracle-Sol\'e, J.~Ruiz, and S.~Shlosman,
  \emph{{Interfaces in the Potts model I: Pirogov--Sinai theory of the
  Fortuin--Kasteleyn representation}}, Commun. Math. Phys. \textbf{140} (1991),
  81--91.

\bibitem{LMR}
L.~Laanait, A.~Messager, and J.~Ruiz, \emph{{Phase coexistence and surface
  tensions for the Potts model}}, Commun. Math. Phys. \textbf{105} (1986),
  527--545.

\bibitem{LT}
B.~Laslier and F.~B. Toninelli, \emph{Lozenge tilings, {G}lauber dynamics and
  macroscopic shape},  (2013), \url{http://arxiv.org/abs/1310.5844}.

\bibitem{LSW01}
G.~F. Lawler, O.~Schramm, and W.~Werner, \emph{One-arm exponent for {2D}
  critical percolation}, Electron. J. Probab. \textbf{7} (2002), {P}aper 2.

\bibitem{LebML}
J.~L. Lebowitz and A.~Martin-L{\"o}f, \emph{On the uniqueness of the
  equilibrium state for {I}sing spin systems}, Commun. Math. Phys. \textbf{25}
  (1972), 276--282.

\bibitem{Nolin}
P.~Nolin, \emph{Near-critical percolation in two dimensions}, Electron. J.
  Probab. \textbf{13} (2008), 1562--1623.

\bibitem{OnsI}
L.~Onsager, \emph{{Crystal statistics. I. A two-dimensional model with an
  order--disorder transition}}, Phys. Rev. \textbf{65} (1944), 117--149.

\bibitem{Pen74}
R.~Penrose, \emph{The r\^ole of aesthetics in pure and applied mathematical
  research}, Bull. Inst. Math. Appl. \textbf{10} (1974), 266--271.

\bibitem{Pen78}
\bysame, \emph{Pentaplexity}, Eureka \textbf{39} (1978), 16--32, reprinted in
  Math. Intellig. 2 (1979), 32--37.

\bibitem{Perk-AY}
J.~H.~H. Perk and H.~Au-Yang, \emph{{Yang--Baxter} equation}, Encyclopedia of
  Mathematical Physics (J.-P. Fran\c{c}oise, G.~L. Naber, and S.~T. Tsou,
  eds.), vol.~5, Elsevier, 2006, pp.~465--473.

\bibitem{Russo}
L.~Russo, \emph{A note on percolation}, Z. Wahrsch'theorie verw. Geb.
  \textbf{43} (1978), 39--48.

\bibitem{Sch00}
O.~Schramm, \emph{Scaling limits of loop-erased walks and uniform spanning
  trees}, Israel J. Math. \textbf{118} (2000), 221--288.

\bibitem{Sch06}
\bysame, \emph{Conformally invariant scaling limits: an overview and collection
  of open problems}, {Proceedings of the International Congress of
  Mathematicians, Madrid} (M.\ Sanz-Sol\'e {\it et al.}, ed.), vol.~{I},
  European Mathematical Society, Zurich, 2007, pp.~513--544.

\bibitem{Seymour-Welsh}
P.~D. Seymour and D.~J.~A. Welsh, \emph{Percolation probabilities on the square
  lattice}, Ann. Discrete Math. \textbf{3} (1978), 227--245.

\bibitem{Smirnov}
S.~Smirnov, \emph{Critical percolation in the plane: conformal invariance,
  {Cardy's} formula, scaling limits}, C. R. Acad. Sci. Paris Ser. I Math.
  \textbf{333} (2001), 239--244.

\bibitem{Smi07}
\bysame, \emph{{Towards conformal invariance of {2D} lattice models}},
  {Proceedings of the {I}nternational {C}ongress of {M}athematicians, {M}adrid,
  2006} (M.~Sanz-Sol\'e {\it et al.}, ed.), vol.~{II}, European {M}athematical
  {S}ociety, Zurich, 2007, pp.~1421--1452.

\bibitem{Smi10}
\bysame, \emph{Conformal invariance in random cluster models. {I. Holomorphic
  fermions in the I}sing model}, Ann. Math. \textbf{172} (2010), 1435--1467.

\bibitem{Smi-icm}
\bysame, \emph{Discrete complex analysis and probability}, {Proceedings of the
  {I}nternational {C}ongress of {M}athematicians, Hyderabad, 2010} (R.~Bhatia,
  A.~Pal, G.~Rangarajan, V.~Srinivas, and M.~Vanninathan, eds.), vol.~{I},
  Hindustan Book Agency, New Delhi, 2010, pp.~595--621.

\bibitem{Smirnov-Werner}
S.~Smirnov and W.~Werner, \emph{Critical exponents for two-dimensional
  percolation}, Math. Res. Lett. \textbf{8} (2001), 729--744.

\bibitem{Sun11}
N.~Sun, \emph{Conformally invariant scaling limits in planar critical
  percolation}, Probability Surveys \textbf{8} (2011), 155--209.

\bibitem{Sykes_Essam}
M.~F. Sykes and J.~W. Essam, \emph{Some exact critical percolation
  probabilities for site and bond problems in two dimensions}, J. Math. Phys.
  \textbf{5} (1964), 1117--1127.

\bibitem{WW_park_city}
W.~Werner, \emph{Lectures on two-dimensional critical percolation}, Statistical
  Mechanics (S.~Sheffield and T.~Spencer, eds.), vol.~16, IAS--Park City, 2007,
  pp.~297--360.

\bibitem{Werner_SMF}
\bysame, \emph{{Percolation et Mod\`ele d'Ising}}, Cours Specialis\'es,
  vol.~16, Soci\'et\'e Math\'ematique de France, Paris, 2009.

\bibitem{ZSWS}
R.~M. Ziff, C.~R. Scullard, J.~C. Wierman, and M.~R.~A. Sedlock, \emph{The
  critical manifolds of inhomogeneous bond percolation on bow-tie and
  checkerboard lattices}, J. Phys. A \textbf{45} (2012), 494005.

\end{thebibliography}
\bibliographystyle{amsplain}

\end{document}